\documentclass[twoside,12pt]{article}
\usepackage{amsmath,amsthm,amssymb,amscd}
\title {Relative mapping class group of $S^p\times D^q$}
\author{Nikolai A. Krylov\\ ~ \\
International University Bremen\\
School of Engineering \& Science\\
P.O. Box 750 561, 28725 Bremen, Germany\\
E-mail:~n.krylov@iu-bremen.de}

\begin{document}

\date {}

\newtheorem{thm}{Theorem}
\newtheorem{cor}{Corollary}
\newtheorem{lem}{Lemma}
\newtheorem{claim}{Claim}
\newtheorem{dfn}{Definition}
\newtheorem{prop}{Proposition}

\def\G              {{S\pi_p(SO(p))}}
\def\DblM           {{\cal D}M}
\def\DblF           {{\cal D}F}
\def\mcgDM          {\tilde{\pi}_0 {{\rm Diff}}(\DblM)}
\def\RmcgM          {\tilde{\pi}_0 {{\rm Diff}(M,{\rm rel}~\partial)}}
\def\RmcgF          {\tilde{\pi}_0 {{\rm Diff}(F,{\rm rel}~\partial)}}
\def\mcgSpSq        {\tilde{\pi}_0 {{\rm Diff}}(S^p\times S^q)}
\def\mcgSSpSq       {\tilde{\pi}_0 {{\rm SDiff}}(S^p\times S^q)}
\def\RmcgSpDq       {\tilde{\pi}_0 {{\rm Diff}}(S^p\times D^q,{\rm rel}~\partial)}
\def\RmcgSpDp       {\tilde{\pi}_0 {{\rm Diff}}(S^p\times D^p,{\rm rel}~\partial)}
\def\K(V)           {\tilde{\pi}_0V {{\rm Diff}(F,\partial)}}
\def\Rhom           {{{\rm Hom}(H_n(F,\partial F),~\G)}}
\def\mcgM           {\tilde{\pi}_0 {{\rm Diff}(M)}}
\def\kerM           {\tilde{\pi}_0S {{\rm Diff}(M)} }
\def\kerDM          {\tilde{\pi}_0S {{\rm Diff}(\DblM)} }
\def\mcgDF          {\tilde{\pi}_0 {{\rm Diff}}(\DblF)}
\def\Aut            {{\rm Aut}~H_p(\DblF)}
\def\SL             {SL(2,\mathbb Z)}
\def\Gv2            {\Gamma_V(2)}

\def\nat            {\mathbb N}
\def\int            {\mathbb Z}
\def\com            {\mathbb C}

\def\da             {\downarrow}
\def\lra            {\longrightarrow}
\def\ra             {\rightarrow}
\def\hra            {\hookrightarrow}
\def\lmt            {\longmapsto}

\def\vf             {\varphi}
\def\lam            {\lambda}
\def\del            {\delta}
\def\eps            {\epsilon}

\def\tp             {\tilde{\psi }}
\def\tj             {\tilde\jmath}
\def\tf             {\tilde{\vf }}

\maketitle

\parskip=2mm

\begin{abstract}
{Algebraic structure of the group of pseudo-isotopy classes of
diffeomorphisms of the trivial disk bundle over the standard
sphere which restrict to the identity map on the boundary is
determined.}
\end{abstract}


\section{Introduction}

Let $M$ be an oriented smooth manifold with non-empty boundary
$\partial M$. We consider the group of pseudo-isotopy classes of
diffeomorphisms of $M$ which restrict to the identity map on the
boundary. Recall that two diffeomorphisms $f_0,f_1\in {\rm
Diff}(M)$ which keep $\partial M$ pointwise fixed are called {\it
pseudo-isotopic (rel boundary)} if there exists a diffeomorphism
$\Phi: M\times I\lra M\times I$ such that $\Phi|_{M\times
\{0\}}=f_0,~\Phi|_{M\times \{1\}}=f_1$ and $\Phi|_{\partial
M\times I}=Id$. We write $f_0\sim f_1$ to indicate that $f_0$ is
pseudo-isotopic to $f_1$. Such diffeomorphisms are of course
orientation preserving and if dim$(M)=2$, the group is known as
the classical mapping class group of a surface. We will also call
the group of pseudo-isotopy classes of diffeomorphisms of $M$
which are fixed on the boundary as {\it the (relative) mapping
class group} and denote it by $\tilde{\pi}_0 {\rm Diff}(M,{\rm
rel}~\partial)$. Not much is known about these groups in higher
dimensions and the goal of this work was to determine such a group
for the trivial disk bundles over the standard spheres (see
Theorem 1 at the end). For the handlebodies in general, the
pseudo-isotopy classes of diffeomorphisms with no constraint on
the boundary had been studied by Wall in \cite{Wall1}. The
corresponding group is denoted by $\tilde{\pi}_0 {\rm Diff}(M)$.
Our approach here will be based on the results of Levine
\cite{Levine1} and Sato \cite{Sato} who determined the mapping
class group of $S^p\times S^q$ up to extension (cf. also work of
Turner, \cite{Turner}).

Integer coefficients are understood for all homology groups,
unless otherwise stated, and symbols $\simeq$ and $\cong$ are used
to denote diffeomorphism and isomorphism respectively. We will
follow notations of \cite{Sato} and denote by $\tilde{\pi}_0 {{\rm
SDiff}}(S^p\times S^q)$ the subgroup of $\tilde{\pi}_0 {{\rm
Diff}}(S^p\times S^q)$ which consists of classes with
representatives that induce trivial action on the homology. The
$h$-cobordism classes of all homotopy $m$-spheres form an abelian
group under the connected sum operation and we denote such a group
by $\Theta_m$ (see \cite{KerMil} for details).

\section{$\tilde{\pi}_0 {\rm Diff}(S^p\times D^q,{\rm rel}~\partial)$}

Given a manifold $M$ with non-empty boundary, one can consider the
{\it double} $\DblM$ of $M$ defined as $\DblM:=\partial(M\times
I)$. $\DblM$ is a closed manifold with the canonically defined
smooth structure (cf. \cite{Munkres}). Since $\partial(M\times
I)\simeq M\times\{0\}\cup (\partial M)\times I\cup M\times\{1\}$
and $(\partial M)\times I\cup M\times\{1\}\simeq M$ (which we will
denote by $M_+$), one can also think of the double as of the union
of two copies of $M$ glued together along the boundary:
$$
\DblM\simeq M\cup M_+
$$
For example, if $M\simeq S^p\times D^q$ then $\DblM\simeq
S^p\times S^q$. Take $\vf\in {\rm Diff}(M,{\rm rel}~\partial)$,
then one can use the identity map to extend $\vf$ to a
diffeomorphism $\tilde{\vf}$ of $\DblM$. To be more precise, we
define
$$
\tilde{\vf}(x):=\left\{
\begin{array}{rcl}
\vf(x)& {\rm if} & x\in M\\
x~& {\rm if} & x\in M_+\\
\end{array}
\right.
$$
This construction gives a map ${\rm Diff}(M,{\rm
rel}~\partial)\lra {\rm Diff}(\DblM)$, which induces a
homomorphism $\omega: \RmcgM\lra \mcgDM$ defined by
$\omega([\vf]):=[\tf]$. The following proposition generalizes
Theorem 2 of \cite{KK}.

\begin{prop} ~ \\
Homomorphism $\omega: \RmcgM\to \mcgDM$ is injective for all $M$.
\end{prop}
\begin{proof}
If $\vf\in \RmcgM$ and $\tf\sim Id$, then there exists an
extension $\Phi\in {\rm Diff}(M\times I)$ of $\tf\in {\rm
Diff}(\partial(M\times I))$. Since $\Phi$ is equal to $\vf$ on
$M\times \{0\}$ and the identity map on $\partial(M)\times I\cup
M\times \{1\}$, this $\Phi$ is a relative pseudo-isotopy that
connects $\vf$ with $Id$.
\end{proof}

It is easy to see that $\mcgSpSq/\mcgSSpSq\cong\int_2$ when $p\ne
q$ (cf. \cite{Sato}, Theorem I or \cite{Levine1}, \S1.2). Since
the extension diffeomorphism $\tf$ fixes $S^p\times D^q_+ =
M_+\subset \DblM$, it follows from Proposition 2.1 of \cite{Sato}
and proposition above that $\RmcgSpDq$ must be a subgroup of
$\mcgSSpSq$. Suppose that $1\leq p<q$ and $3\leq q$, or $4\leq
p=q$, then there is a homomorphism
$$
B: \mcgSSpSq \lra \pi_p(SO(q+1)).
$$
defined by Sato as follows (cf. \cite{Sato}, \S3): Take a
representative $f$ of a class $[f]\in\mcgSSpSq$ and pick a point
$z\in S^q$. Then $S^p=S^p\times z\subset \DblM$ will represent a
generator of $H_p(\DblM)\cong\int$. Since $f$ acts trivially on
$H_p(\DblM)$, it follows from the Hurewicz theorem and the result
of Haefliger \cite{Haefliger1} that there exists a diffeomorphism
$f'\sim f$ which is the identity on $S^p\times z$. Furthermore, if
we take the tubular neighborhood $S^p\times D^q$ of the sphere
$S^p\times z$, then by the tubular neighborhood theorem we can
assume that $f'$ is isotopic to $f''$ such that $f''(x,y) = (x,
b(f)\cdot y)$, where $(x,y)\in S^p\times D^q$ and $b(f)$ is a
smooth map $S^p\lra SO(q)$. Denote the inclusion map $SO(q)\hra
SO(q+1)$ by $S$ and the homotopy class of $b(f)$ by
$[b(f)]\in\pi_p(SO(q))$, then $B$ is defined by the formula:
$$
B([f]):=S_*([b(f)]).$$

For each element $[f]\in\tilde{\pi}_0 {{\rm Diff}}(D^m,{\rm
rel}~\partial)\cong \tilde{\pi}_0 {{\rm Diff}}(S^m)\cong
\Theta_{m+1}$ one defines a diffeomorphism $\iota_r(f)\in {\rm
Diff}(S^p\times D^{m-p},{\rm rel}~\partial)$ as the identity map
outside an embedded disk ${\mathbb D}^m\hra {\rm Int}(S^p\times
D^{m-p})$ and $f|_{D^m}$ on this ${\mathbb D}^m$ (see \S4 of
\cite{Sato} for the details). It is easily deduced from \S4 of
\cite{Sato} that this construction induces a monomorphism
$\iota_r: \Theta_{p+q+1}\hra \RmcgSpDq$.

Furthermore, let us denote by ${\rm FC}^{p+1}_q$, the group of the
pseudo-isotopy classes of orientation preserving embeddings of
$S^q\times D^{p+1}$ in $S^{q+p+1}$. This group was introduced by
Haefliger and the reader will find all the details in \S5 of
\cite{Haefliger2}. Here we only mention that ${\rm
FC}^{p+1}_q\cong \pi_q(SO(p+1))$, when $q<2p$ (see
\cite{Haefliger2}, Corollaries 5.9 \& 6.6).

\begin{lem}
$$
\RmcgSpDq\cong\left\{
\begin{array}{lcl}
\Theta_{p+q+1}\oplus {\rm FC}^{p+1}_q & {\rm if} & 1<p<q\\
\Theta_{p+q+1}\oplus \pi_q(SO(p+1)) & {\rm if} & 1<q<p\\
\end{array}
\right.
$$
\end{lem}
\begin{proof}
Assume first that $p<q$, then we have the exact sequence (see
\cite{Sato}, Theorem II or \cite{Levine1}, Theorem 2.4):
$$
0\lra {\rm FC}^{p+1}_q\oplus \Theta_{p+q+1} \lra \mcgSSpSq
\stackrel{B}{\lra} \pi_p(SO(q+1))\lra 0.
$$

Since the diffeomorphism $\tf$ fixes $S^p\times D^q_+$, it is
clear from definition of $B$ that $B([\tf])=\{0\}$, i.e.
$\RmcgSpDq\subset {\rm Ker}(B)$. Sato had shown (see \cite{Sato},
Lemma 3.3) that for any element $[u]\in{\rm Ker}(B)$ one can find
a representative $u\in {\rm Diff}(S^p\times S^q)$ such that
$u|_{S^p\times D^q_+}=Id$, and therefore $\RmcgSpDq\cong {\rm
Ker}(B)$.

Let now $p$ be larger than $q$. In this case ${\rm
Im}(B)=\pi_q(SO(p+1))$ and for every class $[z]\in \pi_q(SO(p+1))$
we can choose a smooth representative $$r: (D^q, S^{q-1})\lra
(SO(p+1), Id)
$$
and define a relative diffeomorphism $\vartheta$ of $S^p\times
D^q$ by the formula $\vartheta(x,y):=(r(y)\cdot x,y)$, where
$(x,y)\in S^p\times D^q$. Let us also use $\omega$ to denote the
inclusion $\RmcgM\hra \kerDM$. Proof of Proposition 3.2 of
\cite{Sato} shows that the composition $B\circ\omega$ is a
surjection. If we take an element $[\vf_0]\in\RmcgM$ such that
$B([\tf_0])=\{1\}$, then propositions 5.2 and 5.3 of \cite{Sato}
imply that $\tf_0$ (modulo some element of $\Theta_{p+q+1}$ if
needed) can be extended to a diffeomorphism of $S^p\times
D^{q+1}$, i.e. $[\tf_0]\in\Theta_{p+q+1}$ and hence the exact
sequence above implies that $\RmcgSpDq\cong
\pi_q(SO(p+1))\oplus\Theta_{p+q+1}$.
\end{proof}

\begin{lem}
$$
\RmcgSpDq\cong\left\{
\begin{array}{lcl}
\{1\} & {\rm if} & p=1, q=2\\
\Theta_{q+1} \oplus \Theta_{q+2} & {\rm if} & p=1, q\geq 3\\
\Theta_{p+2} \oplus \int_2 & {\rm if} & q=1, p\geq 2\\
\end{array}
\right.
$$
\end{lem}
\begin{proof}
$\tilde{\pi}_0 {{\rm Diff}}(S^1\times D^2,{\rm rel}~\partial) \hra
\tilde{\pi}_0 {{\rm Diff}}(S^1\times S^2)$ by our Proposition 1.
Since $\tf\in {\rm Diff}(S^1\times S^2)$ acts trivially on the
homology, it follows from Theorem 5.1 of \cite{Gluck} that $\tf$
is either pseudo-isotopic to the identity or the diffeomorphism
$T$ of $S^1\times S^2$ defined by $T(t,x):=(t,f(t)\circ x)$ where
$f: S^1\lra SO(3)$ is a smooth generator of
$\pi_1(SO(3))\cong\int_2$. If we had $\tf\sim T$, then $\tf$ would
extend to a diffeomorphism of $S^1\times D^3$, and therefore
would be pseudo-isotopic to the identity map. Thus $\tilde{\pi}_0
{{\rm Diff}}(S^1\times D^2,{\rm rel}~\partial)\cong\{1\}$.

Consider now $\tilde{\pi}_0 {\rm Diff}(S^1\times D^q,{\rm
rel}~\partial)$ with $q\geq3$. Here we also can assume that
$\tilde{\pi}_0 {{\rm Diff}}(S^1\times D^q,{\rm
rel}~\partial)\subset {\rm Ker}(B)$ and according to Proposition
6.3 of \cite{Sato}, the latter group is isomorphic to
$\Theta_{q+1}\oplus \Theta_{q+2}$. For an element $[g]\in
\tilde{\pi}_0 {\rm Diff}(D^q,{\rm rel}~\partial)\cong
\tilde{\pi}_0 {\rm Diff}(S^q)\cong \Theta_{q+1}$, we associate the
diffeomorphism $G \in {\rm Diff}(S^1\times D^q)$ defined by the
formula: $G(x,y):=(x,g(y))$. Thus we obtain a map $K:
\Theta_{q+1}\to\tilde{\pi}_0 {\rm Diff}(S^1\times D^q,{\rm
rel}~\partial)$ which is a monomorphism (see Proposition 6.2 of
\cite{Sato}). Since $\Theta_{q+2}$ is also a subgroup of $\RmcgM$,
we see that $\tilde{\pi}_0 {\rm Diff}(S^1\times D^q,{\rm
rel}~\partial)\cong\Theta_{q+1}\oplus \Theta_{q+2}$.

Gluck had shown in \cite{Gluck} (see \S8 - \S13) that
$\tilde{\pi}_0 {{\rm Diff}}(S^2\times D^1,{\rm
rel}~\partial)\cong\int_2$. As the generator of this group, one
can take the homeomorphism $T$ of $S^2\times D^1$ defined by
$T(x,t):=(f(t)\circ x,t)$, where $[f]\in\pi_1(SO(3))$ is as above.
Since $\Theta_4\cong\{0\}$, we can assume for the rest of our
proof that $p\geq 3$. Since $\RmcgM\hra \kerDM$ and using once
again the generalized Dehn twist $(x,y)\to (\alpha(y)\circ x,y)\in
S^p\times I$ with $[\alpha]=$ a smooth generator of
$\pi_1(SO(p+1))$, it is easy to see that $B\circ\omega$ is an
epimorphism, i.e. ${\rm Im}(B\circ \omega)=\int_2$. Take $\vf\in
{\rm Diff}(S^p\times D^1,{\rm rel}~\partial)$ such that $[\tf]\in
{\rm Ker}(B)$. According to Sato (\cite{Sato}, \S6), we have the
exact sequence
$$
0\lra \Theta_{p+1} \stackrel{K}{\lra} {\rm Ker}(B)
\stackrel{C}{\lra} \Theta_{p+2} \lra 0$$ where map $K$ was just
defined in the paragraph above and $C$ is the inverse map to the
monomorphism $\iota_r$ which has been mentioned at the beginning
(all the details regarding the homomorphism $C$ can be found in
\S4 of \cite{Sato}). If we have $[\tf]\in {\rm Ker}(C)$, then we
can assume that there exists $f\in {\rm Diff}(S^p\times S^1)$ such
that $f\sim\tf$ and $f(x,y)=(g(x),y)$ with $[g]\in \tilde{\pi}_0
{\rm Diff}(D^q,{\rm rel}~\partial)\cong \tilde{\pi}_0 {\rm
Diff}(S^q)\cong \Theta_{q+1}$. In this case we could extend $\tf$
to a diffeomorphism of $S^p\times D^2$, that is $\tf\sim Id$.
Hence $\tilde{\pi}_0 {{\rm Diff}}(S^p\times D^1,{\rm
rel}~\partial)\cong\Theta_{p+2}\oplus \int_2$ as required.
\end{proof}

Consider now a parallelizable $2p-$manifold $F$, which is obtained
by gluing $\mu$ handles of index $p\geq 2$ to the $2p-$disk and
rounding the corners:
$$
F=D^{2p}\cup\bigsqcup_{i=1}^{\mu} (D_i^p\times D^p)$$ Evidently
$S^p\times D^p$ is an example of such a manifold. Given now a
diffeomorphism $\vf$ of $F$ which is the identity map on $\partial
F$, one can consider the variation homomorphism ${\rm Var}(\vf):
H_p(F,\partial F)\lra H_p(F)$ defined by the formula ${\rm
Var}(\vf) [z]:=[f(z)-z]$ for any relative cycle $z\in
H_p(F,\partial F)$ (cf. \S1 of \cite{Stevens}). It is easy to show
(\cite{KK}, \S2.2) that the elements of $\RmcgF$ that induce zero
variation homomorphism (i.e. $[\vf]$ such that ${\rm
Var}(\vf)[z]=0,~\forall [z]$) form a subgroup of $\RmcgF$. We will
follow \cite{KK} and denote this subgroup by $\K(V) $. Let us also
denote by $h$ the homomorphism $\mcgDF \lra \Aut$ induced by the
natural action of the elements of $\mcgDF $ on the $p$-th homology
of the double.

\begin{claim} ~\\
Kernel of the homomorphism $h\circ\omega: \RmcgF\lra \Aut$ is
equal to $\K(V) $.
\end{claim}
\begin{proof}
It follows immediately from the proof of Theorem 1 of \cite{KK}.
\end{proof}

\noindent Let us denote by $\G$ the image of $\pi_p(SO(p))$ under
the map $S_*:\pi_p(SO(p))\lra \pi_p(SO(p+1))$ induced by the
inclusion $SO(p)\hra SO(p+1)$. Then $S\pi_6(SO(6))$ is trivial and
for all other $p\geq 3$ the groups $S\pi_p(SO(p))$ are given in
the table below (see \cite{Kreck}, p. 644):

\parskip=5mm

\begin{tabular}{|c|c|c|c|c|c|c|c|c|}
\hline $p$ (mod 8) & 0 & 1 & 2 & 3 & 4 & 5 & 6 & 7\\
\hline $S\pi_p(SO(p))$ & ~$\int_2\oplus\int_2$~ & ~$\int_2$~ &
~$\int_2$ ~ & ~$\int$ ~ & ~ $\int_2$~ & ~ 0 ~ & ~ $\int_2$ ~ & ~$\int$ \\
\hline
\end{tabular}

~

\begin{lem}
$$
\RmcgSpDp\cong\left\{
\begin{array}{lcl}
\int & {\rm if} & p=1\\
\{1\} & {\rm if} & p=2\\
\Theta_{2p+1}\oplus \G & {\rm if} & 4\leq p~{\rm is~even}\\
\Theta_{2p+1}\oplus \G \oplus \int & {\rm if} & 3\leq p~{\rm
is~odd}
\end{array}
\right.
$$
\end{lem}
\begin{proof}
The case of $p=1$ is well known and a proof can be found, for
example, in \S7 of \cite{Gluck}. When $p$ is even, the image of
$h:\mcgDF\lra{\rm Aut}~H_p(S^p\times S^p)$ is isomorphic to
$\int_4$ and for any element $[\psi]\in\mcgDF$, $h([\psi])$ leaves
no non-zero cycle of $H_p(S^p\times S^p)$ invariant (see
Proposition 2.2 of \cite{Sato}). Since the extended diffeomorphism
$\tf$ preserves $S^p\times D^p_+$ pointwise, it is clear that the
group $\RmcgF$ coincides with the kernel of $h\circ\omega$ and
hence, by the claim above, $\RmcgF\cong\tilde{\pi}_0V {{\rm
Diff}(S^p\times D^p,\partial)}$. The statement for $p=$even will
now follow from the exact sequence (see \cite{Sato}, Theorem II)
$$
0\lra S\pi_p(SO(p)) \oplus \Theta_{2p+1} \lra \tilde{\pi}_0 {{\rm
SDiff}}(S^p\times S^p) \stackrel{B}{\lra} S\pi_p(SO(p))\lra 0,
$$ a simple observation which has been already made that for
every element $[\vf]\in\tilde{\pi}_0V {{\rm Diff}(S^p\times
D^p,\partial)}$, the element $[\tf]$ belongs to the kernel of $B$
and Theorem 3 of \cite{KK} which says that

\noindent If $n=2$ then $\K(V) = 0$, and for all $n\geq 3$ the
following sequence is exact
$$
0\lra \Theta_{2p+1} \stackrel{\iota_r}{\lra} \tilde{\pi}_0V {{\rm
Diff}(F,\partial)} \lra {\rm Hom}(H_p(F,\partial F),~\G) \lra 0.
$$

Assume now that $p$ is odd. Then ${\rm Aut}~H_p(S^p\times
S^p)\cong\SL $ when $p$ is 1, 3 or 7, and in the other cases ${\rm
Aut}~H_p(S^p\times S^p)$ is a proper subgroup of $\SL $ which
consists of the matrices $
\begin{pmatrix}
d_1 & d_2\\
d_3 & d_4\\
\end{pmatrix}$ such that both products $d_1d_2$ and $d_3d_4$ are even
integers (\cite{Wall}, Lemma 5). Clearly, any matrix of this type
is congruent modulo 2 either to $ Id=
\begin{pmatrix}
1 & 0\\
0 & 1\\
\end{pmatrix}$ or
$ V :=\begin{pmatrix}
0 & -1\\
1 & 0\\
\end{pmatrix}$.  We will denote this subgroup by
$\Gv2 $. It is well known that $\Gv2 $ is not a normal subgroup of
index 3 of $\SL $ (see \cite{Rankin}, \S1.5). Moreover, using the
fact that the corresponding projective group
$\Gamma_V(2)/\int_2\cong\int_2 * \int$ is generated by $V$ and
$T:=\begin{pmatrix}
1 & 2\\
0 & 1\\
\end{pmatrix}$
(cf. \cite{Wall}, \S 3) one can easily show that $\Gv2 $ admits
the following presentation $ \Gv2 \cong \langle
V,T~|~V^4=id,~V^2T=TV^2\rangle $. Assume that $p\geq 3$. It
follows again from the definition of $\omega$ that the image of
$h\circ\omega$ consists of those automorphisms that preserve the
class of $H_p(S^p\times S^p)$ represented by an embedded sphere
$S^p\times\{*\}\subset S^p\times D^p_+\subset S^p\times S^p$. If
we choose two spheres $S^p\times\{*\}$ and $\{*\}\times S^p$ as
the basis of $H_p(S^p\times S^p)$, we see that ${\rm
Im}(h\circ\omega)$ is generated either by $\begin{pmatrix}
1 & 1\\
0 & 1\\
\end{pmatrix}$ (when $p=3$ or $p=7$) or by
$\begin{pmatrix}
1 & 2\\
0 & 1\\
\end{pmatrix}$ (in the other cases) and hence ${\rm Im}
(h\circ\omega)\cong\int$. As for the corresponding element of
$\RmcgM$, one can again take the generalized twist $\vartheta$ of
$S^p\times D^p$ defined by the formula
$\vartheta(x,y):=(\zeta(y)\circ x,y)$, where $(x,y)\in S^p\times
D^p$ and $\zeta: (D^p, S^{p-1})\lra (SO(p+1), Id)$ is a smooth map
which generates image of the map $j_*: \pi_p(SO(p+1))\to
\pi_p(S^p)$ from the exact homotopy sequence of the fibration
$SO(p)\hra SO(p+1)\stackrel{j}{\lra} S^p$. To finish the proof we
need to show that for $p=$odd we also have $\tilde{\pi}_0V {{\rm
Diff}(S^p\times D^p,\partial)}\cong \Theta_{2p+1}\oplus \G$. If
$p\geq 5$, one can use exactly the same argument which we gave
above for $p=$even and if $p=3$, it is shown in Example 1 of
\cite{KK} that $\tilde{\pi}_0V {{\rm Diff}(S^3\times
D^3,\partial)}\cong \Theta_7\oplus\int$.
\end{proof}

\noindent Let us now summarize what we have proved above and state
the main result of this paper.
\begin{thm}
$$
\RmcgSpDq\cong\left\{
\begin{array}{lcl}
\int & {\rm if} & p=q=1\\
\{1\} & {\rm if} & 1\leq p\leq 2, q=2\\
\Theta_{p+q+1}\oplus \G & {\rm if} & 4\leq p=q~{\rm is~even}\\
\Theta_{p+q+1}\oplus \G \oplus \int & {\rm if} & 3\leq p=q~{\rm
is~odd}\\
\Theta_{p+q+1} \oplus \Theta_{q+1} & {\rm if} & p=1, q\geq 3\\
\Theta_{p+q+1}\oplus {\rm FC}^{p+1}_q & {\rm if} & 1<p<q\\
\Theta_{p+q+1}\oplus \pi_q(SO(p+1)) & {\rm if} & 1\leq q<p
\end{array}
\right.
$$
\end{thm}

\end{document}